\documentclass[11pt]{article}

\usepackage{amsmath, amsthm, amsfonts, amssymb}
\newtheorem{thm}{Theorem}[section]
\newtheorem{cor}[thm]{Corollary}
\newtheorem{lem}[thm]{Lemma}
\newtheorem{prop}[thm]{Proposition}

\theoremstyle{definition}
\newtheorem{defn}{Definition}[section]

\numberwithin{equation}{section}
\theoremstyle{remarks}
\newtheorem{rem}{Remark}[section]

\title{{\bf Moderate deviations for non-linear functionals and empirical spectral density of moving average processes }}

\author{ H. Djellout\footnote{Laboratoire de Math\'ematiques Appliqu\'ees,
CNRS-UMR 6620, Universit\'e Blaise Pascal, 63177 Aubi\`ere, France,  djellout@math.univ-bpclermont.fr}, A. Guillin\footnote{Ceremade, CNRS - UMR 7534, Universit\'e Paris IX Dauphine, 75775 Paris,France, guillin@ceremade.dauphine.fr}, L. Wu\footnote{Laboratoire de Math\'ematiques Appliqu\'ees, CNRS-UMR 6620, Universit\'e Blaise Pascal, 63177 Aubi\`ere, France,  Li-Ming.Wu@math.univ-bpclermont.fr\ and Department of Mathematics, Wuhan University, 430072 China}}

\date{First Version: Avril 2004}

\newcommand{\ee}{\mathbb{E}}

\newcommand{\nn}{\mathbb{N}}
\newcommand{\rr}{\mathbb{R}}
\newcommand{\pp}{\mathbb{P}}

\newcommand{\ttt}{\mathbb{T}}
\newcommand{\zz}{\mathbb{Z}}

\def\LL{\mathcal L}

\def\<{\langle}
\def\>{\rangle}

\def\beq{\begin{equation}}
\def\neq{\end{equation}}

\def\bdef{\begin{defn}}
\def\ndef{\end{defn}}

\def\bthm{\begin{thm}}
\def\nthm{\end{thm}}

\def\bprop{\begin{prop}}
\def\nprop{\end{prop}}

\def\brmk{\begin{rem}}
\def\nrmk{\end{rem}}

\def\bexa{\begin{exa}}
\def\nexa{\end{exa}}

\def\bexe{\begin{exe}}
\def\nexe{\end{exe}}

\def\bprf{\begin{proof}}
\def\nprf{\end{proof}}

\def\bdes{\begin{description}}
\def\ndes{\end{description}}

\begin{document}

\maketitle

\begin{abstract}
A moderate deviation principle for functionals, with at most quadratic growth, of moving average processes is established. The main assumptions on the moving average process are a Logarithmic Sobolev inequality for the driving random variables and the continuity, or weaker, of the spectral density of the moving average process. We also obtain the moderate deviations for the empirical spectral density, exhibiting an interesting new form of the rate function, i.e. with a correction term compared to the Gaussian rate functionnal.
\end{abstract}
\vspace{20pt}
{\bf AMS 2000 Subject Classification:} 60F10; 60G10; 60G15.\par\vspace{10pt}
{\bf Key Words:}  moderate deviations; moving average processes; logarithmic Sobolev inequalities, toeplitz matrices.\par\vspace{20pt}

\section{Introduction}
Consider the moving average process
\beq
X_n:= \sum_{j=-\infty}^{+\infty} a_{j-n} \xi_j=  \sum_{j=-\infty}^{+\infty} a_j \xi_{n+j},\ \forall n\in\zz.
\label{linear}
\neq

where the innovations $(\xi_n)_{n\in\zz}$ is a sequence of $\rr^d$-valued  centered square integrable i.i.d.r.v., with common law $\LL(\xi_0)=\mu$, and $(a_n)_{n\in\zz}$ be a sequence of real numbers such that
\beq
\sum_{n\in\zz}|a_n|^2<+\infty.
\label{cond}
\neq

This last condition (\ref{cond}) is necessary and sufficient for the a.s. convergence  or convergence in law of the serie (\ref{linear}). The sequence $(X_k)$ is strictly stationary having spectral density
$$f(\theta):={\rm Var}(\xi_0)|g(\theta)|^2$$
where
\beq
g(\theta):=\sum_{n=-\infty}^{+\infty} a_ne^{in\theta}.
\label{dens}
\neq
\par

The moving average processes are of special importance in time series
analysis and they arise in a wide variety of contexts. Applications to
economics, engineering and physical sciences are very broad and a vast
amount of literature is devoted to the study of the limit theorems for
moving average processes under various conditions (e.g. Brockwell and
Davis \cite{BrocDav} and references therein). For example, the {\it minimal} condition for the central limit theorem for $(X_n)$ is (see \cite[Corollary 5.2, p.135]{HH}) that $g$ is continuous at $\theta=0$. The large deviations theorems have attracted much attention and many work, see Burton and Dehling \cite{BuD}, Jiang, Rao and Wang \cite{JRW1},\cite{JRW2}, Djellout and Guillin \cite{DG} and recently by Wu \cite{Wu99} on the linear case, under different assumptions on the law $\xi_0$, and the spectral density function of $X$, see Wu \cite{Wu99}, for relevant reference and more details.
\par\vspace{5pt}

The main purpose of this paper consists to investigate the Moderate
Deviation Principle (in short MDP) for the so-called empirical periodogram of order $n$ of the process $(X_k)$ defined by
\beq
{\cal I}_n(\theta):=\frac 1n \left|\sum_{k=1}^n X_k e^{ik\theta}
\right|^2
\label{per} 
\neq
which are random elements in the space  $L^p(\ttt,d\theta)$ of $p$-integrable function on the torus $\ttt$ identified with $[-\pi,\pi[$ equipped with the weak convergence topology. We present a simple proof under some conditions such as the $L^{q}(\ttt,d\theta)$-boundedness of the spectral density of $(X_k)$ and a Logarithmic Sobolev Inequality (in short LSI) for $\mu$.
\vskip 5pt

The quantity (\ref{per}) is one of the main tools in the study of nonparametric statistical estimation of the unknown spectral density $f$ on the basis of the sample $(X_1,\cdots,X_n)$ from the process $(X_n)$. There exists an abundant literature on several properties and limit theorems of (\ref{per}), specially in Gaussian case. The central limit theorem was generalised by L. Giraitis and D. Surgailis (\cite{GS}) to non Gaussian case and they proved that  $\sqrt{n}({\cal I}_n(h)-\ee {\cal I}_n(h))$ converge in law to normal distribution ${\cal N}(0,\sigma^2)$. In Gaussian case this result was already proved by Avram \cite{Avr} and Fox and Taqqu \cite{FT}.

We also establish the MDP for  additive non-linear functionals of the moving average processes : 
\beq 
\frac 1n \sum_{k=1}^n F(X_k,...,X_{k+l})
\label{target}
\neq 

where $F$ takes its value in $\rr^m$, under some  regularity for the derivatives of $F$. This regularity enables us in particular to obtain the MDP for
$$F(X_k,...,X_{k+l})=\left(X_kX_k^*,X_k X_{k+1}^*,...,X_k X_{k+l}^*\right)$$
which is of particular interest in statistics.
\par\vspace{5pt}

To our knowledge, it is the first time a MDP for functionals of
moving average is established, for a general class of measurable functions
$F$ (and not only in the Gaussian case). Bryc and Dembo \cite{BrD2} have
considered quadratic functional of Gaussian processes both at the level of
large and moderate deviations. We extend their results for the MDP as our
r.v. are not necessarily Gaussian (under the same hypothesis on the
density), and we consider the autocorrelation vector (in a non i.i.d.
setting). Moreover, and compare with Bercu and al \cite{BGR}, we also establish the
MDP for the empirical spectral density,  not only for marginals of the empirical spectral measures. We exhibit an interesting new form of the rate function, i.e. with a correction term compared to the Gaussian rate functional.
\vskip 5pt

Recall that any real stationary Gaussian process $(X_n)$ with a square integrable spectral density function $f$ can be represented as (\ref{linear}), so that one may see our
results as the moderate deviations alternative to the seminal work of
Donsker and Varadhan \cite{DV} on large deviations of Gaussian processes.
\vskip 5pt

\vspace{5pt}
This paper is structured as follows. The MDP for the empirical spectral density is stated in next section. The MDP for non-linear functionals is given is section 3. We establish the key a priori estimation in section 4. The last section is devoted to the proofs of the main results.

\vspace{10pt}


\section{MDP for the empirical spectral density}

In this section we only consider, without loss of generality, and to simplify notations, the real case. Let $(\xi_n)_{n\in\zz}$ is a sequence of $\rr$-valued  centered i.i.d.r.v., with common law $\LL(\xi_0)=\mu$, and let $a:=(a_n)_{n\in\zz}$ be a sequence of real, and define $(X_n)$ by (\ref{linear}). We will always assume that $\mu$ satisfies a LSI, i.e. there exists $C>0$ such that
\beq
{\rm Ent}_\mu(h^2)\le 2C \ee_\mu\left(|\nabla h|^2\right)
\label{logsob}
\neq
for every smooth $h$ such that  $\ee_\mu(h^2\log^+h^2)<\infty$, where
$${\rm Ent}_\mu(h^2)=\ee_\mu(h^2\log h^2)-\ee_\mu(h^2)\log\ee_\mu(h^2).$$
See Ledoux \cite{L99} for further details on LSI. Note that it implies in particular that there exists some positive $\delta$ such that
\beq
\ee_\mu\left(e^{\delta |x|^2}\right)<\infty.
\label{integ}
\neq
Let $(b_n)$ a sequence of real number such that 
\beq 1\ll b_n \ll\sqrt{n} \label{speed}. \neq
For any  measure $\lambda$ on the torus $\ttt$ (identified with $[-\pi,\pi[$, in the usual way), let
$$L^p(\ttt,d\lambda):=\left\{h ~\text{measurable}: ||h||_p=\left(\int_{\ttt}|h(\lambda)|^pd\lambda\right)^{1/p}<\infty\right\},\quad 1\le p <\infty,$$
and 
$$L^{\infty}(\ttt,d\lambda):=\bigg\{h ~\text{measurable}: ||h||_{\infty}=\text{esssup}_{\lambda\in \ttt}|h(\lambda)|<\infty\bigg\}.$$
We are interested in the MDP of the empirical spectral density of $(X_n)$ defined by
$$
{\cal I}_n(\theta):=\frac 1n \left|\sum_{k=1}^n X_k e^{ik\theta}
\right|^2 
$$
which are random elements in the space  $L^p(\ttt,d\theta)$ equipped with the weak convergence topology.

\vskip 5pt

We first present here the MDP for the empirical autocorrelation vector which will be our main tool for the MDP of the empirical spectral density, and has its own interest for statistics. Let $\displaystyle\kappa_4=\frac{\ee(\xi^4)-3\ee (\xi^2)^2}{\ee(\xi^2)^2}.$

\bthm
Suppose that $\mu$ satisfies the LSI (\ref{logsob}), that $(a_n)_{n\in\zz}$ satisfies (\ref{cond}). Suppose moreover that the spectral density function  $f$ is in $L^q(\ttt,d\theta)$, where $2<q\le +\infty$ and $\frac{b_n}{{\sqrt n}}n^{1/q}\rightarrow 0$, then $\displaystyle\left(\frac{1}{b_n\sqrt n}\sum_{k=1}^n\big(X_kX_{k+\ell}-\ee X_kX_{k+\ell}\big)\right)_{0\le \ell\le m}$ satisfies the MDP on $\rr^{m+1}$ with speed $b_n^2$ and with the rate function given by
$$
I(z)=\sup_{\lambda\in \rr^{m+1}}\left\{\langle\lambda,z\rangle -{1\over 2}\lambda^*\Sigma^2\lambda\right\};
$$
where $\Sigma^2=(\Sigma^2_{k,\ell})_{0\le k,\ell\le m}$ and
$$
\aligned
\Sigma^2_{k,\ell}&=\frac{1}{2\pi}\int_{\ttt}\left(e^{i(k-\ell)\theta}+e^{i(k+\ell)\theta}\right)f^2(\theta)d\theta+\kappa_4\left(\frac{1}{2\pi}\int_{\ttt}f(\theta)e^{ik\theta}{\rm d}\theta\right)
\left(\frac{1}{2\pi}\int_{\ttt}f(\theta)e^{i\ell\theta}{\rm d}\theta\right).
\endaligned
$$

\nthm

\brmk
\rm
The additional assumption on the normalizer $b_n$ is exactly the one supposed in Bryc-Dembo \cite[Th. 2.3]{BrD2}, but they only consider the case $l=m=0$ in the Gaussian setting. Their large deviations result (namely Prop. 2.5 in \cite{BrD2}) for the empirical autocorrelation is moreover restricted to the i.i.d. case.
\nrmk

\brmk
\rm
First note that there exists some practical criteria ensuring the fact that a measure $\mu$ satisfies some LSI. For example, consider a $C^2$ function $W$ on $\rr^d$ such that $e^{-W}$ is integrable with respect to Lebesgue measure and let
\begin{equation}
d\mu(x)=Z^{-1} e^{-W(x)}dx\end{equation}
and suppose that for some $c$ in $\rr$, $W''(x)\ge c{\rm Id}$ for every $x$ and that for some $\epsilon>0$,
\begin{equation}\int\int e^{(c^-+\epsilon)|x-y|^2}d\mu(x)d\mu(y)<\infty\label{condlogsob}\end{equation}
where $c^-=-\min(c,0)$. Then $\mu$ satisfies (\ref{logsob}) by the criterion of Wang \cite{L99}. Obviously Gaussian variables fulfill this criterion. See Bobkov-G\"otze \cite{BG} for a necessary and suffient condition in the real case, relying on Hardy's inequalities.
\nrmk

The following corollary follows from Theorem 2.1
\begin{cor} Under the assumptions of Theorem 2.1, we have for all $\ell \ge 0$,\\ $\displaystyle\left({1\over \sqrt{n}b_n}\sum_{k=1}^n \left( X_kX_{k+l}-\ee X_kX_{k+l}\right)\right)$ satisfies the MDP on $\rr$ with speed $b_n^2$ and rate function given by
$$I^{\ell}(z)=\frac{1}{2}\frac{z^2}{\frac{1}{2\pi}\int_{\ttt}(1+\cos(2l\theta))f^2(\theta){\rm d}\theta+\kappa_4\left(\frac{1}{2\pi}\int_{\ttt}f(\theta)\cos(\ell\theta){\rm d}\theta\right)^2}.$$

\end{cor}

\brmk
\rm
Now assume that $(\xi_n)$ is a sequence of real i.i.d. normal random variables, so $(X_n)$ is a stationary Gaussian process and inversely any real Gaussian stationary process $(X_n)$ with a square integrable spectral density function $f$ can be represented as (\ref{linear}).\par
In this case, we have $\ee(\xi^4)=3\ee (\xi^2)^2$ and thus $\kappa_4=0$, so we  obtain

$$I^{\ell}(z)=\frac{1}{2}\frac{z^2}{\frac{1}{2\pi}\int_{\ttt}(1+\cos(l\theta))f^2(\theta){\rm d}\theta}.$$
\nrmk

\vskip 5pt

Let us  present now the main result of this paper. From Theorem 2.1 (and its proof) together with the projective limit method, we yield the functional type's MDP below, for

$$
{\cal L}_n(\theta)=\frac{\sqrt n}{b_n}\left({\cal I}_n(\theta)-\ee{\cal I}_n(\theta)\right).
$$

\bthm
Suppose that $\mu$ satisfies the LSI (\ref{logsob}), that $(a_n)_{n\in\zz}$ satisfies (\ref{cond}). Suppose moreover that the spectral density function  $f\in L^q(\ttt,d\theta)$, where $2<q\le +\infty$ and $\frac{b_n}{{\sqrt n}}n^{1/q+1/{p'}}\rightarrow 0$. Let  $\frac 1p+ \frac{1}{p'}=1$ and $\frac{1}{p'}+\frac 1q< \frac 12$,  then $({\cal L}_n)_{n\ge 0}$ satisfies the MDP on $(L^p(\ttt,d\theta), \sigma(L^p(\ttt,d\theta),L^{p'}(\ttt,d\theta)))$ with speed $b^2_n$  with the rate function given for all  even $\eta\in L^p(\ttt,d\theta)$ by

$$I(\eta)=\left\{\aligned
&\frac{1}{2\pi}\int_{\ttt}\frac{\eta^2(\theta)}{4f^2(\theta)}d\theta-\frac{\kappa_4}{2+\kappa_4}\left(\frac{1}{2\pi}\int_{\ttt}\frac{\eta(\theta)}{2f(\theta)}d\theta\right)^2\\
&\text{~if~}\displaystyle  \eta(\theta)d\theta \text{~is absolutely~continuous~w.r.t.~}f(\theta)d\theta
\text{~and~} \displaystyle  \frac{\eta(\theta)}{f(\theta)} \in L^2(\ttt,d\theta) ;\\
&+\infty ,\quad\quad\text{ otherwise.}\endaligned
\right.
$$

\nthm

\brmk
\rm
Now assume $(X_n)$ is a stationary Gaussian process, so we  obtain that $({\cal L}_n)_{n\ge 0}$ satisfies the MDP on $L^p(\ttt,d\theta)$  with speed $b^2_n$  with the rate function given by

$$I(\eta)=\left\{\aligned
&\frac{1}{2\pi}\int_{\ttt}\frac{\eta^2(\theta)}{4f^2(\theta)}d\theta\\
&\text{~if~} \displaystyle  \eta(\theta)d\theta \text{~is absolutely~continuous~w.r.t.~}f(\theta)d\theta
\text{~and~}\displaystyle \frac{\eta(\theta)}{f(\theta)} \in L^2(\ttt,d\theta)  ;\\
&+\infty ,\quad\quad\text{ otherwise.}\endaligned
\right.
$$

We thus give the MDP for the spectral empirical measure in the setting of Bercu and al \cite{BGR}, note however that they only consider the marginal LDP, i.e. LDP for ${\cal I}_n(h)$ for some bounded $h$ on the torus with an extra assumption on the eigenvalues of the Toeplitz matrix, where
$\displaystyle{\cal I}_n(h)=\frac{1}{2\pi }\int_{\ttt}{\cal I}_n(\theta)h(\theta)d\theta.$
\nrmk

\brmk
\rm
Notice that the extra term with respect to the Gaussian case in the evaluation of the rate function was also found by L. Giraitis and D. Surgailis (\cite{GS}) in their investigations of the CLT for ${\cal I}_n(h)$. The result of (\cite{GS}) can be summarized as below : if
\beq
\lim_{n\rightarrow \infty}\frac 1n {~\rm tr~}((T_n(f)T_n(h))^2)=\frac{1}{2\pi}\int_{\ttt}f^2(\theta)h^2(\theta)d\theta;
\label{GS}
\neq
(where $T_n(h)$ is the Toeplitz matrix of $h$)then $\sqrt{n}({\cal I}_n(h)-\ee {\cal I}_n(h))$ converges in law (as $n\rightarrow \infty$) to the normal distribution ${\cal N}(0,\sigma^2)$ with 
$\displaystyle \sigma^2:=\frac{2}{2\pi}\int_{\ttt}\left(f(\theta)h(\theta)\right)^2d\theta+\kappa_4\left(\frac{1}{2\pi}\int_{\ttt}f(\theta)h(\theta)d\theta\right)^2$.
In Gaussian case this result was already proved by Avram \cite{Avr} and Fox and Taqqu \cite{FT}.
\nrmk

\brmk
\rm
Our main tool in the proof of our Theorem 2.3 is (\ref{GS}), which is valid under our conditions on $f$ and $h$. It seems that the single condition that the integral on the right hand side of (\ref{GS}) is finite (i.e. $h\in L^2(\ttt,f^2d\theta)$) is not  sufficient to obtain (\ref{GS}). This explains why we cannot obtain the MDP of the empirical spectral density in $L^2(\ttt,f^2d\theta)$.

\nrmk

\brmk
\rm One can not hope that the MDP in Theorem 2.3 holds w.r.t. the strong topology of $L^p(\ttt,d\theta)$, because the rate function $I(\eta)$ is not inf-compact w.r.t. this topology.

\nrmk

As a consequence of Theorem 2.3 we have the following

\begin{cor} Under the assumptions of Theorem 2.3, we have  that for all $h\in L^{p'}(\ttt,d\theta)$
$$\limsup_{n\rightarrow \infty}\frac{1}{b_n^2}\log\left(e^{b_n^2\frac{1}{2\pi}\int_{\ttt}h(\theta){\cal L}_n(\theta)d\theta}\right)=\frac 12\left(\frac{2}{2\pi}\int_{\ttt}h^2(\theta)f^2(\theta)d\theta+\kappa_4\left(\frac{1}{2\pi}\int_{\ttt}h(\theta)f(\theta)d\theta\right)^2\right).
$$
\end{cor}

In the next corollary of Theorem 2.3, we replace $\ee{\cal I}_n(\theta)$ by $f(\theta)$, more useful in practice.

\begin{cor} Under the assumptions of Theorem 2.3, assume moreover that $f'\in L^2(\ttt,d\theta)$. The same conclusion holds for $\tilde{\cal L}_n$ instead of ${\cal L}_n$ where

$$
\tilde{\cal L}_n(\theta)=\frac{\sqrt n}{b_n}\left({\cal I}_n(\theta)-f(\theta)\right).
$$
\end{cor}

\brmk\rm
By looking carefully at the proof of this corollary, one may see that the needed convergence of $\ee{\cal I}_n(h)$ to $\int fh$ is ensured by our assumption on $f'$ wich is surely too strong (as the negligibility of this term is in ${1\over \sqrt{n}b_n}$) but remains practical, solely relying on the spectral density. Other possibilities impose implicit, and thus difficult to check, conditions linking $h$ and $f$.    
\nrmk
\section{MDP for non-linear functionals}

Let us present now the following sligthly more general model: $(\xi_n)_{n\in\zz}$ is a sequence of $\rr^d$-valued  centered i.i.d.r.v., with common law $\LL(\xi_0)=\mu$, and let $a:=(a_n)_{n\in\zz}$ be a sequence of real $p\times d$-matrix. We now present the MDP for a functional $F:(\rr^p)^{l+1}\to\rr^m$, i.e. the MDP of $$S_n(F)={1\over \sqrt{n}b_n}\sum_{k=1}^n \left( F(X_k,...,X_{k+l})-\ee\left(F(X_k,...,X_{k+l})\right)\right),$$
and we use the notation $F(x_0,...,x_l)$, so that $\partial_{x_i}F$ should be understood as usual. Let $f(\theta)=g(\theta)\Gamma(\xi_0)g^*(\theta)$, $\Gamma(\xi_0):=({\rm cov}(\xi^i_0,\xi^j_0)_{i,j=1\cdot\cdot\cdot,d}).$

\bthm
Suppose that $\mu$ satisfies the LSI (\ref{logsob}), that $(a_n)_{n\in\zz}$ satisfies (\ref{cond}) and $g$ is continuous on $\ttt$. Suppose moreover that $\partial_{x_i} F$ is Lipschitz for $i=0,...,l$, then $S_n(F)$ satisfies the MDP with speed $b^2_n$ and good rate function $I_F$ given by
$$
I_F(z)=\sup_{\lambda\in \rr^m}\left\{\langle\lambda,z\rangle -{1\over 2}\lambda^*\Sigma^2_F\lambda\right\}=\frac 12 z^*\Sigma_F^{-2}z.
$$
where $\Sigma_F^{-2}$ is the generalized inverse of the covariance matrix $\Sigma_F^{2}$ given by
\beq
\Sigma_F^{2}:=\lim_{n\rightarrow +\infty}\frac 1n \Gamma\left(\sum_{k=1}^nF(X_k,...,X_{k+l})\right)
\label{cc}
\neq
which exists.

\nthm

\brmk
\rm
Note also that under our assumption on $F$ it enables us to obtain the MDP for
$$F(X_k,...,X_{k+l})=\left(X_kX_k^*,X_k X_{k+1}^*,...,X_k X_{k+l}^*\right)$$
as the derivatives in each coordinate is Lipschitz, without further assumption on the normalizer $b_n$ but with a bounded spectral density.
\nrmk

Note also the following corollary in the linear case $F(x_0,..,x_l)=x_0$ which weakens the assumptions on $g$.

\begin{cor}
Suppose that $\mu$ satisfies the integrability condition (\ref{integ}),
that $(a_n)_{n\in\zz}$ satisfies (\ref{cond}) and $g$ is continuous on a
neighborhood of 0, then $S_n(F)$ satisfies a MDP with speed $b^2_n$ and
rate $\displaystyle I(z)=\sup_{\lambda\in \rr^m}\{\langle
z,\lambda\rangle-\frac 12\lambda^* f(0) \lambda\}$.
\end{cor}

It generalizes Th. 3.1 of Djellout and Guillin \cite{DG} to the case of
unbounded r.v. Under assumption (\ref{integ}), the crucial inequality
(\ref{crucial}), as a consequence of the LSI, may not be used. However, we may encompass this difficulty by noting that integrability (\ref{integ}) is, by Djellout and al. \cite[Th. 2.3]{DGW02}, equivalent to a Transport inequality in $L_1$-Wasserstein distance which is itself equivalent to the inequality (\ref{crucial}) with the Lipschitz norm instead of the gradient in the right hand side, but for this particular linear case, the gradient and Lipschitz norm are equal so that the same proof works. The release of the assumptions of the
continuity of $g$ comes from the fact that in this case, Lemma 4.3 is not
used.


\section{A priori estimation}

We recall the following well known elementary result
\begin{lem}
Suppose $Y=[Y_1,\cdot\cdot\cdot,Y_n]^*$ is a real valued centered Gaussian vector with covariance matrix $R$ and let  $A$ be a symmetric real valued $n\times n$-matrix. Then with $\lambda_1,\cdot\cdot\cdot,\lambda_n$ the eigenvalues of the matrix $AR$

\beq
\log \ee \exp(z\langle Y ,AY \rangle)=\left\{\aligned
&-\frac 12\sum_{j=1}^n\log(1-2z\lambda_j)\quad\text{ if }\displaystyle \quad z\max_{1\le j\le n}\lambda_j<1/2\\
&+\infty ,\quad\quad\text{ otherwise.}\endaligned
\right.
\label{gau}
\neq

\end{lem}

We give a crucial lemma which was first stated in Wu \cite{Wu99}, and reproduced here for completeness.

\begin{lem}
If the centered r.v. $\xi_0$ satisfies (\ref{integ}), then there is some constant $K>0$ such that
$$
L(y):=\ee\exp(\langle \xi_0,y\rangle )\le \exp\left(\frac {K^2}2|y|^2\right),\ \forall y\in \rr^d.
$$
\end{lem}

{\bf Proof : } By Chebychev's inequality,
$$
\pp(|\xi_0|>t)\le \exp(-t^2\delta) \ee\exp(\delta|\xi_0|^2):=C(\delta) \exp(-t^2\delta),\ \forall t>0,
$$
consequently
$$
\aligned
L(y)
&\le \ee\exp(|\xi_0||y|)= 1+ \int_0^{\infty} |y|e^{t|y|} \pp(|\xi_0|>t)dt\\
&\le 1+ C(\delta)|y| \int_0^{\infty} \exp(t|y|-t^2\delta) dt\\
&\le 1+ C(\delta) |y| \int_{-\infty}^{\infty} \exp(t|y|-t^2\delta) dt\\
&=  1+ C(\delta)  \sqrt{\frac {\pi}{\delta}}|y| \exp\left(\frac {|y|^2}{4\delta}\right).
\endaligned
$$
Thus there is $C_1>0$ such that (2.1) holds for all $|y|>1$.

For $|y|\le 1$, notice that $\log L(y)\in C^{\infty}(\rr^d)$, and $\log L(0)=0$, $\nabla\log L(y)|_{y=0}=\ee\xi_0=0$. By Taylor's formula of order 2, we have for all $y$ with  $|y|\le 1$,
$$
\log L(y) \le \frac 12 C_2^2 |y|^2,
$$
where $C_2:=\displaystyle\sup_{|y|\le 1} \left(\sum_{k,l=1}^d [\partial_{y_k} \partial_{y_l}\log L(y)]^2\right)^{1/4}$. Thus (2.1)
follows with $K:=C_1\vee C_2$.\hfill$\diamondsuit$

\vspace{10pt}

We extend  (\ref{gau}) from Gaussian distribution to general law $\mu$ satisfying (\ref{integ}), which is a slight generalization of the preceding lemma.  

\begin{lem}
Let $X=[X_1,\cdots,X_n]' \in (\rr^p)^n$ with covariance matrix $A=(A_{k,l})_{1\le k,l \le n}$ where $A_{k,l}$ is a $p\times p$ matrix given by

$$A_{k,l}:=\ee(X_kX_l^*)=\frac{1}{2\pi}\int_{\ttt}e^{i(k-l)\theta}g(\theta)g(\theta)^*d\theta.$$  
Let $B$ be a symmetric real valued $pn\times pn$-matrix. Assume (\ref{integ}). Let $K>0$ given in lemma 4.2 .Then with $\mu^{pn}_1,\cdot\cdot\cdot,\mu^{pn}_{pn}$ the eigenvalues of the matrix $\sqrt{B}A\sqrt{B}$

$$\log \ee \exp(\lambda\langle X ,BX \rangle)\le \left\{\aligned
&-\frac 12\sum_{j=1}^{pn}\log(1-2K^2\lambda\mu^{pn}_j)\quad\text{ if }\displaystyle \quad \lambda\max_{1\le j\le pn}\mu^{pn}_j<\frac{1}{2K^2}\\
&+\infty ,\quad\quad\text{ otherwise.}\endaligned
\right.
$$

\end{lem}

{\bf Proof : } The main difficulty resides in the nonlinear property of $<x,Bx>$. The trick consists to reduce it to an estimation of linear type in the following way :
$$
\ee \left\{ \exp\left[\frac 12 t^2 <X,BX>
\right]\right\}=\ee \left\{ \exp\left[\frac 12 t^2 \left|\sqrt{B}X\right|^2
\right]\right\}
=\int_{(\rr^p)^n} \ee\left\{ \exp\left[t \langle \sqrt{B} X, Y\rangle
\right]\right\}   \gamma(dY)
$$
where $\gamma$ is the standard Gaussian law $N(0,I)$ on $(\rr^p)^n$.

Since
$$
 \langle \sqrt{B}X, Y\rangle  =  \langle X,  \sqrt{B}Y\rangle =
\sum_{k=1}^n \langle X_k,  (\sqrt{B}Y)_k\rangle 
=\sum_{j\in\zz}\langle \xi_j, \sum_{k=1}^n a_{j-k}^*(\sqrt{B}Y)_k\rangle .
$$
where $(a_{j,k}^*)$ is the hermitian transposition of the matrix $(a_{j,k})$, we get by Lemma 4.2.  and the i.i.d. property of $(\xi_j)$,
$$
 \ee\left\{ \exp\left[t \langle \sqrt{B}X, Y\rangle
\right]\right\}  \le \exp\left[\frac {K^2t^2}2 \sum_{j\in\zz}\left|\sum_{k=1}^n a_{j-k}^*(\sqrt{B}Y)_k\right|^2\right].
$$
Now observe that
$$
\aligned
\sum_{j\in\zz}\left|\sum_{k=1}^n a_{j-k}^*(\sqrt{B}Y)_k\right|^2
&= \sum_{k,l=1}^n \sum_{j\in\zz} \langle a_{j-k}^*(\sqrt{B}Y)_k, a_{j-l}^*(\sqrt{B}Y)_l\rangle \\
&=\sum_{k,l=1}^n  \langle (\sqrt{B}Y)_k, \sum_{j\in\zz}a_{j-k}a_{j-l}^*(\sqrt{B}Y)_l\rangle \\
&=\sum_{k,l=1}^n  \langle (\sqrt{B}Y)_k, A_{k,l}(\sqrt{B}Y)_l\rangle \\
&=\langle Y, \sqrt{B} A \sqrt{B}Y\rangle \\
\endaligned
$$

Then letting $\mu^{pn}_1,\cdot\cdot\cdot,\mu^{pn}_n$ be the eigenvalues of the matrix $\sqrt{B}A\sqrt{B}$ (which are also the eigenvalues of $AB$), we get for all $\lambda$ such that $\displaystyle  t^2K^2\max_{1\le j\le pn}\mu^{pn}_j  <1$
$$
\aligned
\ee \left\{ \exp\left[\frac 12 t^2 <X,BX>
\right]\right\}
&\le
\int_{(\rr^p)^n}
\left\{ \exp\left[ \frac 12 K^2t^2 \langle Y, \sqrt{B} A \sqrt{B}Y\rangle 
\right]\right\}
 \gamma(dY)\\
&=-\frac 12\sum_{j=1}^{pn}\log(1-K^2t^2\mu^{pn}_j) 
\endaligned
$$
and it follows with $\lambda =t^2/2$. \hfill$\diamondsuit$

\brmk If we assume $\|g\|_{\infty}=\||g(\theta)|\|_{L^{\infty}(\rr, d\theta)}$, and $B=I$ we obtain exactly the result in Wu \cite{Wu99}. In fact in this case, we have for any $\lambda>0$ such that $2\lambda K^2 \|g\|_{\infty}^2<1$,
\beq
\log \ee  e^{\lambda \langle X,X\rangle}\le -\frac 12 \log \left(1-2\lambda K^2 \|g\|_{\infty}^2 \right)^{np}.
\label{borne}
\neq
\nrmk

\brmk Instead of lemma 4.2., we can use the consequence of the LSI (\ref{crucial}) to prove Lemma 4.3., but (\ref{crucial}) is more stronger than (\ref{integ}) (see below).
\nrmk


\section{Proofs}

Introduce first the following coefficients for each $N\in\nn^*$ :
$\displaystyle a_j^N= a_j \left(1-{|j|\over N}\right)~{\rm if}~|j|\le N~{\rm and}~0~~{\rm otherwise}$, and define the Fejer approximation of $X_k$ and $g$
$$\displaystyle X^N_k=\sum_{j\in \zz}a_j^N \xi_{k+j},\quad g^N(\theta)=\sum_{j\in \zz}a_j^N e^{ij\theta}{~~~\forall \theta \in \rr},$$
that will enable us to first consider the finite case and then extend it to the infinite case by approximation. Remark that if $f\in L^q(\ttt,d\theta)$, $q>2$, then $\int_\ttt (g-g^N)^4d\theta\to0$ as $N\to\infty$.

For any real and symmetric function $h\in L^1(\ttt,d\theta)$, let $T_n(h)$ be the Toeplitz matrix of  $n$ associated with $h$ i.e. $T_n(h)=(\hat{r}_{k-l}(h))_{1\le k,l\le n}$ where $\hat{r}_k(h)$ is the $k$th Fourier coefficient of $h$
$$\hat{r}_k(h)=\frac{1}{2\pi}\int_{\ttt}e^{ik\theta} h(\theta)d\theta, \qquad \forall k \in \zz.$$

The matrix $T_n(h)$ is obviously  real and symmetric, is positive definite whenever $h\ge 0$.
For an $n\times n$ matrix $A$, we consider the usual operator norm $\displaystyle ||A||=\sup_{x\in\rr^n}\frac{|Ax|}{|x|}$.

\vspace{5pt}

We shall need the two following lemmas. The first gives an estimate for the maximal eigenvalue of the covariance matrices $T_n(f)$ which is Lemma 4.7 of Bryc-Dembo \cite{BrD2}. The second one concerning the asypmtotic behavior of the trace of the products of Toeplitz matrices see (\cite{GS}). \par

\begin{lem} If $1\le q \le \infty$ then  for all $n>1$ we have $||T_n(f)||\le  n^{1/q}||f||_q.$
\end{lem}

\begin{lem} Let $f_k\in L^1(\ttt,d\theta)\cap L^{q_k}(\ttt,d\theta)$ with $0\le q_k \le \infty$ for $k=1,\cdots s$ and $\displaystyle\sum_{k=1}^s \frac{1}{q_k}\le 1$. The following assertion hold
$$\lim_{n\rightarrow \infty}\frac{1}{n}{\rm~tr~}\left(\prod_{k=1}^sT_n(f_k)\right)=\hat{r}_0\left(\prod_{k=1}^s f_k\right).$$ 
\end{lem}

\subsection{Proof of Theorem 2.1}

We shall prove it only in the real valued case. The proof is divided into three steps. In the first one, we prove that the MDP holds for some suitable approximation of our process, then we will show this approximation is a good one in the sense of the moderate deviations and we will finally establish the convergence of the rate function and the subsequent existence of the limiting variance.\par\vspace{10pt}

{\it Step 1}.
Let
$$Q_n^N={1\over \sqrt{n}b_n}\sum_{k=1}^n \left(X^N_kX^N_{k+l}-\ee X^N_kX^N_{k+l})\right).$$
The crucial remark is that the sequence $X^N_kX^N_{k+l}$ is a $2N$-dependent identically distributed sequence. Using (\ref{integ}), we get for all $N$ and for some positive $\eta$ that $\ee\left(e^{\eta |X^N_kX^N_{k+l}|}\right)<\infty.$

We may then apply results of Chen \cite{Ch} on Banach valued MDP of $m-$dependent sequence, enabling us to get that for each $N$ fixed, for all $\lambda$

\begin{equation}
\begin{array}{l}
\displaystyle \lim_{n\rightarrow \infty}\frac{1}{b_n^2}\log\ee  \left(e^{\lambda b_n^2 Q_n^N } \right)=\displaystyle\frac{\lambda^2}{2}\lim_{n\rightarrow\infty}\frac {1}{n}\ee \left(\sum_{k=1}^n X^N_kX^N_{k+l}-\ee X^N_kX^N_{k+l}\right)^2 \\
\quad\qquad\qquad\qquad\qquad\qquad:=\displaystyle\frac{\lambda^2}{2}\Sigma^2_N \in \rr ,
\label{ch2}\\
\quad\qquad\qquad\qquad\qquad\qquad=\displaystyle\frac{\lambda^2}{2}\sum_{k=-N}^{N}{\rm Cov}\left(X^N_0X^N_l, X^N_kX^N_{k+l}\right)
\end{array}
\end{equation}

and that $Q^N_n$ satisfies the MDP with the good rate function
$\displaystyle I^N(x)=\sup_{\lambda\in \rr}\left\{\lambda x -\frac {\lambda^2}{2} \Sigma^2_N\right\}~.$
\vspace{10pt}

{\it Step 2}. The purpose of this step will be to prove the asymptotic negligibility as $N\to\infty$ of $Q_n-Q^N_n$ with respect to the MDP, i.e. we will establish that for all $\lambda\in \rr$
\beq
\lim_{N\to\infty}\limsup_{n\to\infty}{1\over b^2_n}\log \ee\left(e^{\lambda b_n^2(Q_n-Q^N_n)}\right)=0.
\label{negli2}
\neq
Remark that, by Jensen inequality and as our functionals are centered, we only have to establish the upper inequality in (\ref{negli2}).

Our main tool is the following consequence of the LSI (\ref{logsob}), see Ledoux \cite[Th. 2.7]{L99} applied to our context (after having extended (\ref{logsob}) by tensorization to the infinite product measure of $\mu$): for exponentially integrable $G$,
\beq
\ee\left( e^{\lambda {b_n\over \sqrt{n}}(G-\ee G)}\right)\le \ee \left( e^{\lambda^2 {b_n^2\over n}C|\nabla G|^2}\right),
\label{crucial}
\neq
with $C$ given in (\ref{logsob}). Let apply it to $$G((\xi_l)_{l\in\zz})=\sum_{k=1}^n(X_kX_{k+l}-X^N_kX^N_{k+l}),$$ so that our main estimations are now transferred to the gradient of $G$.

Clearly

$$\partial_{\xi_i}G=\sum_{k=1}^n(a_{i-k}X_{k+l}+a_{i-k-l}X_k-a_{i-k}^NX_{k+l}^N-a_{i-k-l}^NX_k^N) ;$$

so

$$
\aligned
|\nabla G|^2 &\le 4 \sum_{i\in \zz}\left((\sum_{k=1}^n(a_{i-k}-a^N_{i-k})X_{k+l})^2+(\sum_{k=1}^n(a_{i-k-l}-a^N_{i-k-l})X_k)^2\right.\\
&+\left.(\sum_{k=1}^na^N_{i-k}(X_{k+l}-X^N_{k+l}))^2+(\sum_{k=1}^na^N_{i-k-l}(X_k-X^N_k))^2\right)\\
&=(I)+(II)+(III)+(IV).
\endaligned
$$

By H\"older inequality,

\beq
\aligned
\log\ee\left(e^{\lambda \frac{b_n} {\sqrt n} (G-\ee G)}\right) &\le \log\ee\left(e^{C\lambda^2 \frac{b^2_n} {n} ||\nabla G||^2}\right) \\
&\le \frac 14 \log\ee\left(e^{4C\lambda^2 \frac{b^2_n} {n} (I)}\right)+\frac 14 \log\ee\left(e^{4C\lambda^2 \frac{b^2_n} {n} (II)}\right)\\
&+ \frac 14 \log\ee\left(e^{4C\lambda^2 \frac{b^2_n} {n} (III)}\right)+\frac 14 \log\ee\left(e^{4C\lambda^2 \frac{b^2_n} {n} (IV)}\right).\\
\label{es}
\endaligned
\neq

Let us deal with the first term of this inequality. We rewrite the expression of $(I)$ as

$$
\aligned
(I) &= 4 \sum_{i\in \zz}\sum_{k,k'=1}^n(a_{i-k}-a^N_{i-k})(a_{i-k'}-a^N_{i-k'})X_{k+l}X_{k'+l}\\
&=4 \sum_{k,k'=1}^n \hat{r}_{k'-k}((g-g^N)^2) X_{k+l}X_{k'+l}\\
&=4<X_{\cdot+l},T_n((g-g^N)^2)X_{\cdot+l}>
\endaligned
$$

Let $\mu_1^{n,N},\cdot\cdot\cdot,\mu_n^{n,N}$ be the eigenvalues of the matrix $$\sqrt{T_n((g-g^N)^2)}T_n(f)\sqrt{T_n((g-g^N)^2)}.$$
Its operator norm is bounded from above by (using Lemma 5.1) 
$$||T_n(f)||\cdot||T_n((g-g^N)^2)||\le n^{1/q}||f||_q n^{1/q}||(g-g^N)^2||_q.$$
Since $\frac{b_n}{\sqrt n} n^{1/q}\to 0$ and $f\in L^q(\ttt,d\theta)$, we choose $n$ sufficiently large such that  $\displaystyle 32C^2\lambda^2\frac{b_n^2}{n} \max_{1\le j\le n}\mu^{n,N}_j<1$. Applying Lemma 4.3, we get

\beq
\log\ee(e^{4C\lambda^2\frac{b_n^2}{n}(I)}) \le -\frac{1}{2}\sum_{j=1}^n\log(1-32CK^2\lambda^2\frac{b_n^2}{n}\mu^{n,N}_j).
\label{es1}
\neq

Similarly, we have

\beq
\aligned
\log\ee(e^{4C\lambda^2\frac{b_n^2}{n}(II)})&= \log\ee e^{16C\frac{b_n^2}{n}\lambda^2<X{\cdot},T_n((g-g^N)^2) X_{\cdot}>}\\
&\le -\frac{1}{2}\sum_{j=1}^n\log(1-32CK^2\lambda^2\frac{b_n^2}{n}\mu_j^{n,N}).
\endaligned
\label{es2}
\neq

Let us deal with the third term. We rewrite the expression of $(III)$ as

$$
\aligned
(III) &= 4 \sum_{i\in \zz}\sum_{k,k'=1}^na^N_{i-k}a^N_{i-k'}(X_{k+l}-X^N_{k+l})(X_{k'+l}-X^N_{k'+l})\\
&=4 \sum_{k,k'=1}^n \hat{r}_{k'-k}((g^N)^2)(X_{k+l}-X^N_{k+l})(X_{k'+l}-X^N_{k'+l})\\
&=4<X_{\cdot+l}-X^N_{\cdot+l},T_n((g^N)^2)(X_{\cdot+l}-X^N_{\cdot+l})>.
\endaligned
$$

Let $\nu_1^{n,N},\cdot\cdot\cdot,\nu_n^{n,N}$ the eigenvalues of the matrix $$\sqrt{T_n((g^N)^2)}T_n((g-g^N)^2)\sqrt{T_n((g^N)^2)}.$$
Its operator norm is bounded from above by (using Lemma 5.1) 
$$||T_n((g^N)^2)||\cdot||T_n((g-g^N)^2)||\le n^{1/q}||(g^N)^2||_q n^{1/q}||(g-g^N)^2||_q.$$ 
By our assumptions on $b_n$ and $f$, once again we take $n$ sufficiently large such that  $\displaystyle 32CK^2\lambda^2\frac{b_n^2}{n} \max_{1\le j\le n}\nu^{n,N}_j <1$. Applying lemma 4.3., we get

\beq
\log\ee(e^{4C\lambda^2\frac{b_n^2}{n}(III)})\le-\frac{1}{2}\sum_{j=1}^n\log(1-32CK^2\lambda^2\frac{b_n^2}{n}\nu^{n,N}_j).
\label{es3}
\neq

Similarly

\beq
\aligned
\log\ee(e^{4C\lambda^2\frac{b_n^2}{n}(IV)})\le -\frac{1}{2}\sum_{j=1}^n\log(1-32CK^2\lambda^2\frac{b_n^2}{n}\nu^{n,N}_j).
\endaligned
\label{es4}
\neq

By (\ref{es}) and the previous estimations (\ref{es1}) (\ref{es2}) (\ref{es3}) (\ref{es4}), we have

\beq
\aligned
\frac{1}{b_n^2}\log\ee\left( e^{\lambda b_n^2(Q_n-Q_n^N)}\right)
\le -\frac 14\frac 1n\sum_{j=1}^n\left(\log(1-32CK^2\lambda^2\frac{b_n^2}{n}\mu^{n,N}_j)
+\log(1-32 CK^2\lambda^2\frac{b_n^2}{n}\nu^{n,N}_j)\right).
\label{cr}
\endaligned
\neq

Notice that by the Taylor's expansion of order 1, we have for $|z|<1$
$$\log(1-z)=-z(1-tz)^{-1}$$
where $t=t(z)\in [0,1]$. This applied here to $z^{n,N}_j=32CK^2\lambda^2\frac{b_n^2}{n}\lambda^{n,N}_j$, ($\lambda^{n,N}_j=\nu^{n,N}_j$ or $\lambda^{n,N}_j=\mu^{n,N}_j$) which satisfies $\displaystyle\sup_{1\le j\le n}|z^{n,N}_j|\rightarrow 0$ as $n \rightarrow \infty$, and hence $|1-t(z^{n,N}_j)z^{n,N}_j|\rightarrow 1$ uniformly in $1\le j\le n$. Thus

$$
\lim_{n\rightarrow \infty}\frac{1}{b_n^2}\log \ee\left( e^{\lambda b_n^2 (Q_n-Q_n^N)}\right)\le 16C^2\lambda^2\lim_{n\rightarrow \infty}\left(\frac 1n\sum_{j=1}^n(\mu^{n,N}_j+\nu^{n,N}_j)\right).$$

Thanks to the elementary formula ${\rm ~tr~}(AC)={\rm ~tr~}(CA)$ and using Lemma 5.2, we have
$$
\lim_{n\rightarrow \infty}\frac 1n \sum_{j=1}^n\mu^{n,N}_j=\lim_{n\rightarrow \infty}\frac 1n {\rm ~tr~}\left( T_n(f)T_n((g-g^N)^2)\right)= \hat{r}_0\left((g-g^N)^2f\right).
$$

Similarly

$$
\lim_{n\rightarrow \infty}\frac 1n \sum_{j=1}^n\nu^{n,N}_j=\lim_{n\rightarrow \infty}\frac 1n {\rm ~tr~}\left(T_n((g^N)^2)T_n((g-g^N)^2)\right)=  \hat{r}_0\left((g^N)^2(g-g^N)^4\right).
$$

So we have
$$
\aligned
&\limsup_{n\rightarrow \infty}\frac{1}{b_n^2}\log \ee\left( e^{\lambda b_n^2(Q_n-Q_n^N)}\right)\le 32C^2\lambda^2 \hat{r}_0(f^2)\hat{r}_0((g-g^N)^4).
\endaligned
$$
Letting $N$ to infinity , we get the desired negligibility (\ref{negli2}).

 We then obtain that $Q_n$ satisfies the MDP of speed $b^2_n$ and good rate function $\widetilde I$ by the approximation lemma \cite[Th. 2.1]{Wu99}, with $\widetilde I$ given by

\begin{equation}
\widetilde I(x)=
\displaystyle
\sup_{\delta>0}\liminf_{N\rightarrow \infty}\inf_{B(x,\delta)}I^N
=\sup_{\delta>0}\limsup_{N\rightarrow \infty}\inf_{B(x,\delta)}I^N~.
\label{t2}
\end{equation}

{\it Step 3}.
We have now to prove the identification of the rate function. First, we show that
\begin{equation}
\Sigma^2:=\displaystyle\lim_{n\rightarrow \infty}\frac{1}{n}\ee \left(\sum_{k=1}^n(X_kX_{k+l}-\ee X_kX_{k+l})\right)^2 {\rm ~exists~ and}\quad
\Sigma^2=\lim_{N\rightarrow +\infty}\Sigma_N^2 \in \rr .
\label{cv2}
\end{equation}

By the previous estimations, we have that for all $|\lambda|$ small enough
$$
\aligned
\ee \left(e^{\lambda(G-\ee G)}\right)\le 1
&+16C^2\lambda^2 n \hat{r}_0(f^2)\hat{r}_0((g-g^N)^4)+o\left(\frac{\lambda^2}{2}\right)
\endaligned
$$

Since, for all $|\lambda|$ small enough
$\displaystyle \ee \left(e^{\lambda (G-\ee G) }\right)=1+\frac{\lambda^2}{2}\ee  \left(G-\ee G\right)^2+o(\frac{\lambda^2}{2}),$

we deduce that $\displaystyle
\ee  \left( G-\ee G \right)^2\leq 16C^2n\hat{r}_0(f^2)\hat{r}_0((g-g^N)^4)+o\left(\frac{\lambda^2}{2}\right)
.$

So we have
$$\sup_n\frac{1}{n}\ee  \left(G-\ee G\right)^2\longrightarrow 0~~{\rm as}~N\rightarrow +\infty~.$$

Whence the limit $\Sigma^2$ in (\ref{cc}) exists, and  $\Sigma_N^2\longrightarrow \Sigma^2$.

\vspace{10pt}

Now we claim that

\begin{equation}
\lim_{n\rightarrow \infty}\frac{1}{b_n^2}\log\ee  \exp\left ( \lambda b_n^2 Q_n\right)=\frac {\lambda^2}{2}\Sigma^2.
\label{lim02}
\end{equation}

For fixed $p,q>1$ with $\frac{1}{p}+\frac{1}{q}=1$, by the H\"older inequality we have that

$$
\log\ee  \exp\left (\lambda b_n^2 Q_n\right)\leq
\frac{1}{q}\log\ee  \exp\left (q\lambda b_n^2(Q_n-Q^N_n)\right)+ \frac{1}{p}\log\ee  \exp\left (p b_n^2 \lambda Q^N_n\right)
$$

for all $\lambda$. From (\ref{ch2}) and previous estimations it follows that for some constant $B>0$
$$
\limsup_{n\rightarrow \infty}\frac{1}{b_n^2}\log\ee \left( e^{b_n^2\lambda Q_n}\right)\leq \frac{p\lambda^2}{2}\Sigma_N^2+q B\lambda^2\hat{r}_0((g-g^N)^4).$$

Letting $N\rightarrow \infty$ and using (\ref{cv2}), we get

\begin{equation}
\limsup_{n\rightarrow \infty}\frac{1}{b_n^2}\log\ee \left( e^{b_n^2\lambda Q_n}\right)\leq \frac{p\lambda^2}{2}\Sigma^2.
\label {lim12}
\end{equation}

Similarly, by the H\"older inequality, we have

$$
\log\ee  \exp\left (b_n^2 \lambda Q^N_n \right)\leq
\frac{1}{q}\log\ee  \exp\left (\frac{qb_n^2}{p} \lambda (Q_n^N-Q_n)\right)+\frac{1}{p}\log\ee\exp\left (b_n^2\lambda Q_n \right)$$

for every $\lambda$. From (\ref{ch2}) and previous estimations it follows that

$$
\frac{\lambda^2}{2p^2}\Sigma_N^2\leq
\liminf_{n\rightarrow \infty}\frac{1}{pb_n^2}\log\ee \left( e^{b_n^2\lambda Q_n}\right)+\frac{q\lambda^2}{2p^2}B\hat{r}_0((g-g^N)^4).$$

Letting $N\rightarrow \infty$ and using (\ref{cv2}), we obtain

\begin{equation}
\frac{\lambda^2}{2p}\Sigma^2\leq \liminf_{n\rightarrow \infty}\frac{1}{b_n^2}\log\ee \left( e^{b_n^2\lambda Q_n}\right).
\label{lim22}
\end{equation}

Letting $p\rightarrow 1$ in (\ref{lim12}) and (\ref{lim22}) yields (\ref{lim02}).\par
\vspace{10pt}

So by (\ref{lim02}) and the Laplace principle \cite[Th. 2.1.10, p.43]{DS}, we have

\begin{eqnarray}
\displaystyle \lim_{n\rightarrow \infty}\frac{n}{b_n^2}\log\ee \left( e^{b^2_n\lambda Q_n}\right)&=&\frac {\lambda^2}{2}\Sigma^2\nonumber\\
&=&\sup_{x\in \rr}\big\{ xy -\widetilde I(x)\big\}.\label{lim32}
\end{eqnarray}

To conclude, we have now to show that $\widetilde I(x)$ defined in (\ref{t2}) is convex.

$$\widetilde I\left(\frac{1}{2}(x_1+x_2)\right)=\sup_{\delta>0}\limsup_{N\rightarrow \infty}\inf_{B(\frac{1}{2}(x_1+x_2),\delta)}I^N$$

\begin{eqnarray*}
\inf_{B(\frac{1}{2}(x_1+x_2),\delta)}I^N&\leq& \inf_{y_1\in B(x_1,\delta), y_2\in B(x_2,\delta)}I^N\left(\frac{1}{2}(y_1+y_2)\right)\\
&\leq&\frac{1}{2} \inf_{y_1\in B(x_1,\delta),y_2\in B(x_2,\delta)}\left(I^N(y_1)+I^N(y_2)\right)\\
&=&\frac{1}{2}\left(\inf_{B(x_1,\delta)}I^N+\inf_{B(x_2,\delta)}I^N\right)
\end{eqnarray*}

So \qquad\qquad $\displaystyle
\limsup_{N\rightarrow \infty}\inf_{B(\frac{1}{2}(x_1+x_2),\delta)}I^N\leq \frac{1}{2}\left(\limsup_{N\rightarrow \infty}\inf_{B(x_1,\delta)}I^N+\limsup_{N\rightarrow \infty}\inf_{B(x_2,\delta)}I^N \right)$

Letting $\delta \downarrow 0$, we get $\widetilde I\left(\frac{1}{2}(x_1+x_2)\right)\leq \frac{1}{2}\bigg(\widetilde I(x_1)+\widetilde I(x_2)\bigg)~.$

Since $\widetilde I$ is inf-compact and convex, by Fenchel's theorem and (\ref{lim32}), we get for all $x \in \rr $
$$\widetilde I(x)=\sup_{\lambda\in \rr}\{\lambda x -\frac {\lambda^2}{2} \Sigma^2\},$$
which is exactly the announced rate function.


\subsection{Proof of Theorem 2.3}

We begin with the following lemma \cite[Chap.2, Prop. 2.5]{Wu97} which implies the exponential tightness.

\begin{lem}
Under the hypothesis Theorem 2.3, we have that for all $h\in L^{p'}(\ttt,d\theta)$
$$\limsup_{n\rightarrow \infty}\frac{1}{b_n^2}\log\left(e^{b_n^2\frac{1}{2\pi}\int_{\ttt}h(\theta){\cal L}_n(\theta)d\theta}\right)<+\infty$$
In particular $\pp({\cal L}_n\in \cdot)$ is exponentially *-tight in $(L^p(\ttt,d\theta),\sigma(L^p(\ttt,d\theta),L^{p'}(\ttt,d\theta)))$, where $\frac{1}{p'}+\frac 1p =1$.
\end{lem}

{\bf Proof :} For every function $h\in L^{p'}(\ttt,d\theta)$, the function $\tilde h(\theta)=\frac 12 [h(\theta)+h(-\theta)]$ is even and
$$\frac{1}{2\pi}\int_{\ttt}h(\theta){\cal I}_n(\theta)d\theta=\frac{1}{2\pi}\int_{\ttt}{\tilde h}(\theta){\cal I}_n(\theta)d\theta,$$
we shall hence restrict oureselves to the case where $h$ is even. Since

$$\frac{1}{2\pi}\int_{\ttt}h(\theta){\cal L}_n(\theta)d\theta=\frac{1}{b_n\sqrt n}\left(\langle X_.,T_n(h)X_.\rangle-\ee\langle X_.,T_n(h)X_.\rangle\right)$$

Let apply (\ref{logsob}) to
$H((\xi_l)_{l\in\zz})=\langle X_.,T_n(h)X_.\rangle$ :

$$\ee(e^{\lambda b_n^2(\frac{1}{2\pi}\int_{\ttt}h(\theta){\cal L}_n(d\theta)})=\ee(e^{\lambda \frac{b_n}{\sqrt n}(H-\ee H)})\le \ee(e^{\lambda^2\frac{b_n^2}{n}C|\nabla H|^2}).$$

Clearly

$$
\aligned
|\nabla H|^2=\sum_{i\in\zz}(\partial_{\xi_i}H)^2&=\sum_{i\in\zz} \left(2\sum_{l,k=1}^n  a_{i-k}X_{l}T_n(h)_{k,l}\right)^2\\
&=4\sum_{l,k,l',k'=1}^n T_n(f)_{k,k'} X_l X_{l'}T_n(h)_{k,l}T_n(h)_{k',l'}\\
&=4 \langle X_.,T_n(h)T_n(f)T_n(h)X_.\rangle.
\endaligned
$$

Let $\alpha^n_1,\cdot\cdot\cdot,\alpha^n_n$ the eigenvalues of the matrix $$\sqrt{T_n(h)T_n(f)T_n(h)}T_n(f)\sqrt{T_n(h)T_n(f)T_n(h)}.$$ 

Its operator norm is bounded from above by (using Lemma 5.1)
$$||T_n(f)||\cdot||T_n(h)T_n(f)T_n(h)||\le (n^{1/q}||f||_q)^2(n^{\frac{1}{p'}}||h||_{p'})^2$$
Since $\frac{b_n}{\sqrt n} n^{1/q+1/{p'}}\to 0$, $f\in L^q(\ttt,d\theta)$ and $h\in L^{p'}(\ttt,d\theta)$, we take $n$ large enough such that  $\displaystyle 8CK^2\lambda^2\frac{b_n^2}{n}\max_{1\le j\le n}\alpha_j^n  <1$. Applying Lemma 4.3. we get

$$\log\ee(e^{\lambda b_n^2(\frac{1}{2\pi}\int_{\ttt}h(\theta){\cal L}_n(d\theta)})\le -\frac{1}{2}\sum_{j=1}^n\log\left(1-8CK^2\lambda^2\frac{b_n^2}{n}\alpha^n_j\right).$$

Thus
$$
\limsup_{n\rightarrow \infty}\frac{1}{b_n^2}\log\left(e^{b_n^2(\frac{1}{2\pi}\int_{\ttt}h(\theta){\cal L}_n(d\theta)}\right)
\le 8C^2\lambda^2\lim_{n\rightarrow \infty}\frac{1}{n}\sum_{j=1}^n\alpha^n_j.
$$

Since $f\in L^q(\ttt,d\theta)$ and $h\in L^{p'}(\ttt,d\theta)$ with $\frac{1}{p'}+\frac  1q<\frac 12$, applying Lemma 5.2, we obtain 
$$\lim_{n\rightarrow +\infty}\frac{1}{n}\sum_{j=1}^n\alpha^n_j =\lim_{n\rightarrow +\infty}\frac 1n {\rm ~tr~}\left((T_n(f)T_n(h))^2\right)=\hat{r}_0(f^2h^2)<+\infty.$$
The proof of the Lemma ends.

\vspace{10pt}

We may now turn to the proof of Theorem 2.3.

{\bf Proof :}\par

{\it Step 1.} Since $\displaystyle \left(\frac{1}{b_n{\sqrt n}}\sum_{k=n-\ell+1}^{n}\big(X_kX_{k+\ell}-\ee X_kX_{k+\ell}\big)\right)_{0\le \ell \le m}$ is negligible with respect to the MDP, using Theorem 2.1, we get the finite dimensional MDP on $\rr^{m+1}$ of
$$\left(\frac{1}{b_n{\sqrt n}}\sum_{k=1}^{n-\ell}\big(X_kX_{k+\ell}-\ee X_kX_{k+\ell}\big)\right)_{0\le \ell \le m}$$
with the rate function given by
$$
I(z)=\sup_{\lambda\in \rr^{m+1}}\left\{\langle\lambda,z\rangle -{1\over 2}\lambda^*\Sigma^2\lambda\right\}.
$$

Now notice that
$$\widehat {\cal L}_n(\ell):=\frac{1}{2\pi}\int_{\ttt}e^{i\ell\theta}{\cal L}_n(d\theta)=\frac{1}{b_n{\sqrt n}}\sum_{k=1}^{n-\ell}\big(X_kX_{k+\ell}-\ee X_kX_{k+\ell}\big).$$
Thus $(\widehat {\cal L}_n(\ell))_{0\le \ell\le m }$ satisfies the MDP on $\rr^{m+1}$ with the same rate function. By Lemma 4.3 and the projective limit Theorem \cite[Th. 4.6.9]{DZ}, we deduce that $({\cal L}_n)_{n\ge 0}$ satisfies the MDP on $(L^p(\ttt,d\theta),\sigma(L^p(\ttt,d\theta),L^{p'}(\ttt,d\theta)))$  with the rate function given by for even function $\eta\in L^p(\ttt,d\theta)$ 

\beq
I(\eta)=\sup_{m\ge 0}\sup_{\lambda_0,..,\lambda_m\in \rr}\left\{\frac{1}{2\pi}\int_{\ttt}\left(\sum_{k=0}^m e^{ik\theta}\lambda_k\right)\eta(\theta) d\theta-\frac 12\Lambda\left(\sum_{k=0}^m e^{ik\theta}\lambda_k\right)\right\}
\label{app}
\neq
where
$$
\aligned
\Lambda\left(\sum_{k=0}^m e^{ik\theta}\lambda_k\right)=\lambda^*\Sigma^2\lambda&=\frac{1}{2\pi}\int_{\ttt}\left(\left|\sum_{k=0}^m e^{ik\theta}\lambda_k\right|^2 +\left(\sum_{k=0}^m e^{ik\theta}\lambda_k\right)^2\right)f^2(\theta)d\theta\\
&\qquad\qquad+\kappa_4\left(\frac{1}{2\pi}\int_{\ttt}\left(\sum_{k=0}^m e^{ik\theta}\lambda_k\right)f(\theta)d\theta\right)^2.
\endaligned
$$

{\it Step 2. Identification of the rate function}. Remark as trigonometric polynomials are dense in $L^2(\ttt,f^2d\theta)$ , one can find for $h\in L^2(\ttt,f^2d\theta)$, an approximation by some trigonometric polynomials sequence $h_n$, such that

\beq\lim_{n\rightarrow \infty}\int_{\ttt}\big(h_n-h\big)^2(\theta)f^2(\theta)d\theta=0.
\label{seq}
\neq

So we can extend continuously the definition of $\Lambda$ to all function $h\in L^2(\ttt,f^2d\theta)$  
$$
\Lambda(h)=\frac{2}{2\pi}\int_{\ttt}h^2(\theta)f^2(\theta)d\theta+\kappa_4\left(\frac{1}{2\pi}\int_{\ttt}h(\theta)f(\theta)d\theta\right)^2.
$$

(a) Suppose that  $\eta(\theta)d\theta$ is absolutely continuous w.r.t. $f^2(\theta)d\theta$, and  $\displaystyle\frac{\eta}{f}\in L^2(\ttt,d\theta)$. Let $h_n$ the sequence defined below in (\ref{seq}), by Cauchy-Schwartz inequality, we get for all even function $\eta\in L^{p}(\ttt,d\theta)$

$$\left(\int_{\ttt}\left|(h_n-h)(\theta)\eta(\theta)\right|d\theta\right)^2\le \int_{\ttt}|h_n(\theta)-h(\theta)|^2f^2(\theta)d\theta\int_{\ttt}\left(\frac{\eta}{f}\right)^2(\theta)d\theta\underset{n\rightarrow \infty} \longrightarrow 0.$$

So  $I(\eta)$ defined in (\ref{app}) coincides with 

$$
I(\eta)=\sup_{h\in L^2(\ttt,f^2d\theta)}\left\{\frac{1}{2\pi}\int_{\ttt}h(\theta) \eta(\theta)d\theta-\frac 12\Lambda(h)\right\}:=\sup_{h\in L^2(\ttt,f^2d\theta) }D(h).
$$

Let us find explicitly the maximizer $h_0$ of $D(h)$. Let $k\in L^2(\ttt,f^2d\theta) $ and $\epsilon >0$,
$$
\aligned
\lim_{\epsilon \rightarrow 0}\frac{D(h+\epsilon k)-D(h)}{\epsilon}
=&
\frac{1}{2\pi}\int_{\ttt}k(\theta)\eta(\theta)d\theta-\frac{1}{2}\left(\frac{2}{2\pi}\int_{\ttt}2f^2(\theta)h(\theta)k(\theta)d\theta\right.
\\
&\qquad \left.+2\kappa_4\left(\frac{1}{2\pi}\int_{\ttt}f(\theta)h(\theta)d\theta\right)\left(\frac{1}{2\pi}\int_{\ttt}f(\theta)k(\theta)d\theta\right)\right)
\endaligned
$$

So 
\beq
\displaystyle\lim_{\epsilon \rightarrow 0}\frac{D(h+\epsilon
k)-D(h)}{\epsilon}=0,~~\forall k \in L^2(\ttt,f^2d\theta) 
\label{derive}
\neq
implies that
\beq
\eta(\theta)=2f(\theta)^2h(\theta)+\kappa_4 \left(\frac{1}{2\pi}\int_{\ttt}f(\theta)h(\theta)d\theta\right)f(\theta).
\label{eta}
\neq

Dividing (\ref{eta}) by $f$ and integrating over $\ttt$ , we obtain

$$\int_{\ttt}f(\theta)h(\theta)d\theta=\frac{1}{2+\kappa_4}\int_{\ttt}\frac{\eta(\theta)}{f(\theta)}d\theta.$$

Replacing this last expression in (\ref{eta}), it is then easy to verify that the only functional $h_0\in L^2(\ttt,f^2d\theta)$ realizing (\ref{derive}) is given by
$$h_0(\theta)f(\theta)=\frac{\eta(\theta)}{2f(\theta)}-\frac{\kappa_4}{2+\kappa_4}\left(\frac{1}{2\pi}\int_{\ttt}\frac{\eta(u)}{2f(u)}du\right).$$
Calculating $D(h_0)$ gives finally the announced rate function.

\vspace{5pt}

(b) Now we have to treat the case where $\eta(\theta)d\theta$ is absolutely continuous w.r.t. $f^2(\theta)d\theta$ but  $\displaystyle\frac{\eta}{f}\not\in L^2(\ttt,d\theta)$. So there exists $g\in L^2(\ttt,d\theta)$ such that $\displaystyle\int_{\ttt}g(\theta)\frac{\eta}{f}(\theta)d\theta=+\infty$, and $\displaystyle g\frac{\eta}{f}\ge 0$. Let $\displaystyle h:=\frac{g}{f},$ so $h\in L^2(\ttt,f^2d\theta) $, we choose $h_n=(h\vee(-n))\wedge n$. We get by dominated convergence
$$\lim_{n\rightarrow \infty}\int_{\ttt}\big(h_n(\theta)-h(\theta)\big)^2f(\theta)^2d\theta=0,$$
so it follows that
$$\lim_{n\rightarrow +\infty}\Lambda(h_n)=\Lambda(h).$$

By Fatou's lemma we get
$$\liminf_{n\rightarrow \infty}\int_{\ttt}h_n(\theta)\eta(\theta) d\theta\ge \int_{\ttt}\liminf_{n\rightarrow \infty}h_n(\theta)\eta(\theta) d\theta=+\infty.$$
Since
$$
I(\eta) \ge \frac{1}{2\pi}\int_{\ttt} h_n(\theta)\eta(\theta)d\theta-\frac{1}{2}\Lambda(h_n),
$$
letting $n$ to $\infty$, we obtain $I(\eta)=\infty.$

\vskip 5pt

(c)  Now we have to treat the case where $\eta(\theta)d\theta$ is not absolutely continuous w.r.t. $f^2(\theta)d\theta$, i.e. there exists a set $K\subset \ttt$ such that $\displaystyle\int_{K}f^2(\theta)d\theta=0$ while $\displaystyle\int_K\eta(\theta)d\theta>0$. For any $t>0$, we approximate the function $t1_K$ by a sequence function $h_n\in L^2(\ttt,f^2d\theta) $. So $\forall t\in \rr$

$$
I(\eta)\ge\lim_{n\rightarrow +\infty}D(h_n)\ge t \int_K\eta(\theta)d\theta.
$$
Letting $t$ to infinity, we get $I(\eta)=+\infty$.
\vskip 5pt

\subsection{Proof of corollary 2.5}

Here we assume $f'\in L^2(\ttt,d\theta)$, so $\displaystyle\sum_{k} |k|^2|\hat {r}_k(f)|^2<\infty.$

We thus only need to prove that for all $h\in L^{p'}(\ttt,d\theta)$ ( so $h\in L^2(\ttt,d\theta)$ since $p'\ge 2$)

$$\frac{\sqrt {n}}{b_n}\left(\int_{\ttt}h(\theta)\ee{\cal I}_n(\theta)d\theta-\int_{\ttt} f(\theta)h(\theta)d\theta\right)\underset{n\rightarrow 0}\longrightarrow 0$$

We have
$$
\aligned
\left|\int_{\ttt}h(\theta)\ee{\cal I}_n(\theta)d\theta-\int_{\ttt} f(\theta)h(\theta)d\theta\right|=&\left|\sum_{|k|\le n-1}\left(1-\frac kn\right)\hat {r}_k(f)\hat{r}_k(h)-\sum_{k}\hat {r}_k(f)\hat{r}_k(h)\right|\\
=&\left|-\sum_{|k|\le n-1}\frac{|k|}{n}\hat {r}_k(f)\hat{r}_k(h)-\sum_{|k|\ge n}\hat {r}_k(f)\hat{r}_k(h)\right|.\\ 
\endaligned
$$

We have $\displaystyle\sum_{|k|\ge n}|\hat {r}_k(f)||\hat{r}_k(h)|\le \sum_{|k|\ge n}\frac {|k|}{n}|\hat {r}_k(f)||\hat{r}_k(h)|$. So applying Cauchy-Schwartz inequality we get
$$\aligned
\left|\int_{\ttt}h(\theta)\ee{\cal I}_n(\theta)d\theta-\int_{\ttt} f(\theta)h(\theta)d\theta\right|&\le \frac 1n  \sqrt{\sum_{k} |k|^2|\hat {r}_k(f)|^2}\sqrt{\sum_{k}|\hat {r}_k(h)|^2}\\
\qquad &\le \frac Cn.
\endaligned
$$
The proof ends.


\subsection{Proof of Theorem 3.1}

For simplicity, we only consider the problem in $\rr$ and $F(x_0,...,x_l)=F(x_0)$.

Let us describe briefly how the preceding proof of Theorem 2.1 can be easily extended to the more general framework of our example.

Since $F'$ is Lipschitz continuous, we get for some positive $L$, and for all $N$
$$
|F(X^N_k)|\le L(1+|X^N_k|^2)\le 2L(N+1)\left(1+\sum_{j=-N}^N a_j^2 \xi_{k+j}^2\right)$$
so that, setting $\delta'={\delta\over 2L(N+1)^2\sup_j a_j^2}$ where $\delta$ is given in (\ref{integ}), by the assumption on the validity of the LSI, we get
$$\ee\left(e^{\delta' |F(X_k^N)|}\right)\le e^{\delta' L(N+1)}~\ee\left(e^{\delta \xi_0^2}\right)<\infty.$$

Since Chen \cite{Ch} deals with moderate deviations of $m-$dependent Banach space valued random variables, so that the first step is exactly the same in the general case.

\vspace{10pt}

To prove the asymptotic negligibility as $N\to\infty$ of $S_n(F)-S_n^N(F)$ with respect to the MDP, we  need to assume the boundedness of the density. We apply again (\ref{crucial}) to 
$$G((\xi_l)_{l\in\zz})=\sum_{k=1}^n(F(X_k)-F(X^N_k)), $$  

We have
\begin{eqnarray*}
|\nabla G|^2
&=&\sum_{i\in\zz}\left(\sum_{k=1}^n a_{i-k}F'(X_k)-a^N_{i-k}F'(X^N_k)\right)^2\\
&\le&2\sum_{i\in\zz}\left(\sum_{k=1}^n (a_{i-k}-a^N_{i-k})F'(X_k)\right)^2+2\sum_{i\in\zz}\left(\sum_{k=1}^n a^N_{i-k}(F'(X^N_k)-F'(X_k))\right)^2\\
&=& 2\left|\sqrt {T_n((g-g^N)^2)}F'(X.)\right|^2+2\left|\sqrt {T_n((g^N)^2)}(F'(X.)-F'(X.^N))\right|^2.
\end{eqnarray*}
By the fact that the derivative of $F$ is Lipschitz and the spectral density is bounded, we have that the last term is bounded by
$$ 2 L \|g-g^N\|^2_\infty (n+\langle X_{\cdot},  X_{\cdot}\rangle )+ 2\|g^N\|^2_\infty \langle X_{\cdot}^N-X_{\cdot},X_{\cdot}^N-X_{\cdot} \rangle .$$
Finally by (\ref{borne}), as $\lambda^2 {b^2_n\over n} \|g-g^N\|^2_\infty$ can be chosen arbitrary small for large $n$, 
$$
\aligned
{1\over b^2_n}\log \ee\left(e^{\lambda b^2_n (S_n(F)-S_n(F^N))}\right)&\le LC\lambda^2 \|g-g^N\|^2_\infty\\
& -{n\over 4b_n^2}\log\left(1-4CLK^2\lambda^2 {b^2_n\over n} \|g-g^N\|^2_\infty \|g\|^2_\infty\right)\\
&-{n\over4b_n^2}\log\left(1-4CLK^2\lambda^2 {b^2_n\over n} \|g^N\|^2_\infty\|g-g^N\|^2_\infty\right)
\endaligned
$$
and the left hand side of this last inequality is easily seen to behave as $n\to\infty$ as

$$\|g-g^N\|^2_\infty\left(LC\lambda^2+2CLK^2\lambda^2\|g\|^2_\infty\right).$$

By the famous Fejer Theorem, under the assumption of continuity of $g$, we get that
$$\lim_{N\to\infty}\|g-g^N\|^2_\infty=0,$$
which yields to the desired negligibility.

A careful reading of {\it Step 3} in the proof of Theorem 2.1 shows that the extension to the general case brings no further difficulties. The proof then ends.

\brmk
\rm
To prove negligibility of {\it Step 2} in general framework, we  only have to establish this negligibility for each of the coordinates $F_j$ of $F$ (as there is only a finite number of coordinates), and also that
\begin{eqnarray*}
|\nabla G|^2&=&
\sum_{i\in\zz}\sum_{j=1}^m\left(\partial_{\xi_i}\sum_{k=1}^nF_j(X_k,\cdot\cdot\cdot,X_{k+l})-F_j(X^N_k,\cdot\cdot\cdot,X^N_{k+l}))\right)^2\\
&=&\sum_{i\in\zz}\sum_{j=1}^m\left(\sum_{s=0}^l\sum_{k=1}^n \left(a_{i-k-s}\partial_{x_s}F_j(X_k,\cdot\cdot\cdot,X_{k+l})-a^N_{i-k-s}\partial_{x_s}F_j(X^N_k,\cdot\cdot\cdot,X^N_{k+l})\right)\right)^2\\
&\le&(l+1)\sum_{s=0}^l\sum_{j=1}^m\sum_{i\in\zz}\left(\sum_{k=1}^n \left(a_{i-k-s}\partial_{x_s}F_j(X_k,\cdot\cdot\cdot,X_{k+l})-a^N_{i-k-s}\partial_{x_s}F_j(X^N_k,\cdot\cdot\cdot,X^N_{k+l})\right)\right)^2
\end{eqnarray*}
which leads to the same estimation as before as $\partial_{x_s}F_j$ is supposed to be Lipschitz for each $j$ and $s$.
\nrmk


\end{document}